\documentclass[doublespacing]{elsart}

\newtheorem{theorem}{Theorem}
\newtheorem{lemma}{Lemma}

\newtheorem{corollary}{Corollary}

\newtheorem{remark}{Remark}
\newtheorem{example}{Example}
\input epsf
\usepackage{amssymb}
\usepackage{amsfonts}
\usepackage{amsmath}
\begin{document}

\begin{frontmatter}



\title{    An explicit formulation for two
 dimensional vector partition functions \thanksref{now}}
\thanks[now]{Project Supported by The
National Natural Science Foundation of China(10401021).}

 \author{Zhiqiang Xu}\footnote{Email: xuzq@lsec.cc.ac.cn}
 \address{Institute of Computational Mathematics, Academy of Mathematics and System Sciences, Chinese
 Academy of Sciences, Beijing, China}


\begin{abstract}
Based on discrete truncated powers, the beautiful Popoviciu's
formulation for restricted integer partition function is
generalized. An explicit formulation for two dimensional
multivariate truncated power functions is presented. Therefore, a
simplified explicit formulation for two dimensional vector
partition functions is given. Moreover, the generalized Frobenius
problem is also discussed.

\end{abstract}

\begin{keyword}
Multivariate truncated powers\sep Vector partition function
\end{keyword}
\end{frontmatter}

% main text
\section*{1. Introduction}

The {\it vector partition function } that we are interested in is
in the form of
\[
t({\bf b}|M)=\#\{{\bf x}\in {\bf Z}_+^n|{\bf Mx}={\bf b}\},
\]
where, ${\bf Z}_+$ denotes the nonnegative integers, ${\bf M}$ is
a fixed $s\times n$ integer matrix with columns $m_1,\cdots,m_n\in
{\bf Z}^s$ and ${\bf b }$ is a variable vector in ${\bf Z}^s.$ To
guarantee $t({\bf b}|M)$ is finite, we require
$[\{m_1,\cdots,m_n\}]$ does not contain the origin, where $[A]$
denotes the convex hull of a given set $A$. The vector partition
function $t({\bf b}|M),$ which is also called a {\it discrete
truncated power}, has many applications in different mathematical
areas including Algebraic Geometry \cite{pommersheim},
Representation Theory\cite{repre}, Number Theory \cite{number} ,
Statistics\cite{Astati} and Randomized Algorithm \cite{rando} etc.
.

When $s=1,$ an explicit formulation for  $t({\bf b}|M)$, which
counts the integer solutions for the linear Diophantine equation,
is presented in \cite{BeckFrobenius}. In particular, when
$M=(a,b)$ where $a$ and $b$ are relatively prime, Popoviciu gave a
beautiful and surprising formulation for $t(n|(a,b))$
(\cite{popoviciu}).

 For the general matrix $M$, the nature of $t({\bf b}|M)$ is investigated and
the piecewise structure of $t({\bf b}|M)$ is given in
\cite{dahmen2} and \cite{sturmfels}. Moreover, one is also
interested in the  explicit formulation of $t({\bf b}|M)$. For the
general matrix $M$, a powerful method for obtaining $t({\bf b}|M)$
is described  in \cite{brion,vergne1}. Another interesting
algorithm for computing $t({\bf b}|M)$ as a function of ${\bf b}$
is also introduced in \cite{dgbeck}.  When $M$ is unimodular, in
which every nonsingular square submatrix has determinant $\pm 1$,
two algebraic algorithms for generating the explicit formulation
for $t({\bf b}|M)$ is presented in \cite{deloera}. But all these
methods depend on the complex computation.  In \cite{xu}, based on
multivariate truncated power functions $T({\bf x}|M)$,  an
explicit formulation for $t({\bf b}|M)$ is presented. But the
formulation involves multivariate truncated power functions
$T({\bf x}|M)$, which is not explicit form, and high-dimensional
Fourier-Dedekind sums, so we have to give an explicit form for
$T({\bf x}|M)$ and simplify high-dimensional Fourier-Dedekind
sums, in order to predigest the explicit formulation for $t({\bf
b}|M).$

The rest of the paper is organized as follows.  To help make this
paper self-contained we shall first introduce some notations and
definitions  in Section 2. In Section 3, we recall some results
about  vector partition functions $t({\bf b} |M)$. Section 4
generalize the Popoviciu's formulation. In Section 5, the
generalized Frobenius problem is investigated.  Finally, Section 6
give an explicit formulation for multivariate truncated powers in
the case where $s=2$ and show the high-dimensional
Fourier-Dedekind sum can be converted to one-dimensional
Fourier-Dedekind sum, which is convenient for computing. And
hence,  a simplified  explicit formulation for two-dimension
vector partition functions is given.

\section*{2.Preliminaries}

 To describe the nature of $t({\bf b}|{\bf M}),$ we
introduce several notations and  definitions in which the common
terminology of multiset theory is adopted. Intuitively, a {\it
multiset} is a set with possible repeated elements; for instance
$\{2,2,2,3,4,4\}$. Let $Y$ be an $s\times n$ matrix. $Y$ can be
considered as a multiset of elements of ${\bf R}^s.$
 The cone spanned by $Y,$ denoted
by $cone(Y)$, is the set
$$
\{ \sum_{y\in Y}a_y y:a_y\geq 0 \mbox{\, for all\, } y \}.
$$
Denote by $cone^\circ (Y)$ the relative interior of $cone(Y).$ Let
${\mathcal Y}(M)$ denote the set consisting of those submultisets
$Y$ of $M$ for which $M \backslash Y$ does not span ${\bf R}^s.$
Let the set $c(M)$ be the union of cone$(M\setminus Y)$ where $Y$
runs over ${\mathcal Y}(M).$ A connected component of $cone^\circ
(M)\setminus c(M),$ is called a {\it fundamental $M-$cone}.
 For the fundamental $M-$cone $\Omega,$ we
set $v(\Omega|M):=\Omega-[[M)).$ Here,  $
[[M)):=\{\sum_{j=1}^na_jm_j:0\leq a_j< 1, \forall j\},
\Omega-[[M))$ is the set of all elements of the form $a-b,$ where
$a\in \Omega$ and $b\in [[M)).$

We shall use the standard multiindex notation. Specifically, an
element $\alpha\in {\bf N}^m$ is called an {\it $m-$index}, and
$|\alpha |$ is called the length of $\alpha.$ Define $ z^\alpha
:=z_1^{\alpha_1}\cdots z_m^{\alpha_m}$ for $z=(z_1,\cdots,z_m)\in
{\bf C}^m$ and $\alpha=(\alpha_1,\cdots,\alpha_m)\in {\bf N}^m.$
 For $y=(y_1,\cdots, y_s) \in {\bf R}^s$ and a function
$f$ defined on ${\bf R}^s$ , we denote by $D_yf$ the directional
derivative of $f$ in the direction $y$,i.e. $
D_y=\sum_{j=1}^sy_jD_j, $ where, $D_j$ denote the partial
derivative with respect to the $j$th coordinate. For
$v:=(v_1,\cdots,v_m)\in {\bf N}^m,$ we let $D^v=D_1^{v_1}\cdots
D_m^{v_m}$ and $v!=\prod _i v_i!$ Moreover, we  let
$e:=(1,1,\cdots,1)\in {\bf Z}^s.$

Let ${\mathcal S}_k(M)=\{ Y\subseteq M:\#Y=s+k, span(Y)={\bf
R}^s\}$ and ${\mathcal B}(Y)=\{ X\subseteq Y:\#X=s, span(X)={\bf
R}^s\}.$ If for any $Y\in {\mathcal S}_k(M),$ $gcd \{
|det(X)|,X\in {\mathcal B}(Y) \}=1$ , then $M$ is called a {\it
$k-$ prime matrix}. In particular, when $M$ is an $1-$prime
matrix, $M$ is also called a {\it pairwise relative prime matrix}.
When $s=1,$ $k-$prime matrix means that no $k$ of the integers
$m_1,m_2,\cdots,m_n$ have a common factor, where
$m_i,i=1,\cdots,n$ are the elements in $M.$ Moreover, we denote
$e^{\frac{2\pi i}{k}}$ by $W_k$.

 The multivariate truncated power
$T(\cdot |M)$ associated with $M,$  which was introduced by
W.Dahmen \cite{dahemenpower} firstly, is the distribution given by
the rule
\begin{equation}\label{powerdefinition}
T(\cdot |M):\phi \mapsto \int _{{\bf R}_+^n}\phi (Mu)du , \phi \in
{\mathcal D}({\bf R}^s),
\end{equation}
where ${\mathcal D}({\bf R}^s)$ is the space of test functions on
${\bf R}^s$,i.e. the space of all compactly supported and
infinitely differentiable functions on ${\bf R}^s$. In fact,
$T(\cdot |M)$ agrees with some homogeneous polynomial of degree
$n-s$  on each fundamental $M-$cone. When $M$ is an $s\times s$
invertible matrix, $T(\cdot |M)$ agrees with the function on ${\bf
R}^s$ which takes value $\frac{1}{|det(M)|}$ on $cone(M)$ and $0$
elsewhere.

In \cite{numeircalT}, an efficient method for computing the
multivariate truncated power is presented.

\begin{theorem}(\cite{numeircalT})\label{Th:numT}
Let $M$ be an $s\times n$  matrix with columns  $m_1,\cdots,m_n\in
{\bf Z}^s\setminus \{ 0\}$  such that the origin does not contain
in $conv(M).$ For any $\lambda_1,\cdots,\lambda_n\in {\bf R}$, and
${\bf x}=\sum_{j=1}^n\lambda_j m_j,$
\begin{equation}\label{eq:recT}
T({\bf x}|M)=\frac{1}{n-s}\sum_{j=1}^n\lambda_jT(x|M \setminus
m_j).
\end{equation}
\end{theorem}
 For more detailed information about the function, the
reader is referred to \cite{deboor2},\cite{dahemenpower}.

A multivariate Box spline $B(\cdot |M)$ associated with $M$ was
introduced in \cite{deboor1} and \cite{deboor2}, which is the
distribution given by the rule
\begin{equation}\label{Eq:boxspline}
 B(\cdot
|M):\phi \mapsto \int _{{[0,1)}^n}\phi (Mu)du , \phi \in {\mathcal
D}({\bf R}^s).
\end{equation}
By taking $\phi=exp(-iy\cdot )$ in (\ref{Eq:boxspline}), we obtain
the Fourier transform of $B(\cdot |M)$ as
$$
\widehat{B}(\zeta
|M)=\prod_{j=1}^n\frac{1-exp(-i\zeta^Tm_j)}{i\zeta^Tm_j},\zeta \in
{\bf C}^s.
$$
For more detail information about Box splines,the reader is
referred to \cite{deboorbook} .

\begin{remark}
The definition of fundamental M-cone is slightly different with
the one presented in \cite{dahmen2}. In \cite{dahmen2}, a
fundamental M-cone is defined as a connected component of
$cone^\circ (M)\setminus c(M),$ where $c(M)$ is the union of
$span(M\setminus Y)$ and $Y$ runs over ${\mathcal Y}(M).$ In fact,
the fundamental $M-$cone defined in this paper may be larger  than
the one defined in \cite{dahmen2}. But all the conclusions in
\cite{dahmen2} hold for the larger fundamental M-cone. In a
private communication, Prof. M. Vergne introduce the new
definition about the fundamental M-cone.
\end{remark}

\section*{3. Vector partition functions}

To describe  the nature of $t({\bf b}|M),$ we let $ M_{\theta}:=\{
y\in M:\theta^y=1\}$ and let $ A(M):=\{ \theta\in ({\bf
C}\setminus \{ 0\})^s:span(M_{\theta})={\bf R}^s \}. $ Recall
$e=(1,1,\cdots,1)\in {\bf Z}^s.$ Obviously, $e\in A(M).$

The following qualitative result about $t(\cdot |M)$ is presented
in \cite{dahmen2}.

\begin{theorem}(\cite{dahmen2})\label{Theorem:DaMi,truncted power's structure}
Let $M=\{m_1,\cdots,m_n\}$ be a multiset of integer vectors in
${\bf R}^s$ such that $M$ spans ${\bf R}^s$ and the convex hull of
$M$ does not contain the origin. Then for any fundamental $M-$cone
$\Omega,$ there exists a unique element $f_{\Omega}(\alpha
|M)=\sum\limits_{\theta\in
A(M)}\theta^{\alpha}p_{\theta,\Omega}(\alpha)$ such that
$f_{\Omega}(\alpha |M)$ agrees with $t(\alpha |M)$ on $v(\Omega
|M),$ where $p_{\theta,\Omega}(\cdot )$ is a polynomial with
degree less than $\# M_{\theta}-s.$
\end{theorem}

An explicit formulation for $p_{e,\Omega}(\alpha)$, which is the
polynomial part of $t(\alpha|M),$ is presented in the following
theorem.
\begin{theorem}(\cite{xu})\label{Th:FormP}
Under the condition of Theorem \ref{Theorem:DaMi,truncted power's
structure}, $p_{e,\Omega}(x)=\sum_{k=0}^{n-s}p_{k,\Omega}(x),$
where $p_{k,\Omega}(x)$ is homogeneous polynomial of degree
$n-s-k,$ defined inductively by
$$p_{0,\Omega}(x)=T({\bf x}|M),
p_{k,\Omega}(x)=-\sum_{j=0}^{k-1}(\sum_{|v|=k-j}D^vp_{j,\Omega}(x)
(-i)^{|v|}D^v\widehat{B}(0|M)/v!),k\geq 1,$$ where, $x\in \Omega
.$
\end{theorem}

More generally, an explicit formulation for $p_{\theta,\Omega}$ is
also given as follows.
\begin{theorem}(\cite{xu})\label{Th:FormPtheta}
Given $\theta_0\in A(M)\setminus e,$ under the condition of
Theorem \ref{Theorem:DaMi,truncted power's structure},
$p_{\theta_0,\Omega}(x)=\sum_{\mu=0}^{n-s-\kappa}p^{\theta_0}_{\mu,\Omega}(x)$,
where $\kappa=\# (M\setminus
M_{\theta_0}),p^{\theta_0}_{\mu,\Omega}(x)$ is homogeneous
polynomial of degree $n-s-\kappa-\mu,$ defined inductively by
$$p^{\theta_0}_{0,\Omega}(x)=q_{0,r}^{\theta_0}(x),
p^{\theta_0}_{\mu,\Omega}(x)=q_{\mu,r}^{\theta_0}(x)
-\sum_{j=0}^{\mu-1}(\sum_{|v|=\mu-j}D^vp^{\theta_0}_{\mu,\Omega}(x)
(-i)^{|v|}D^v\widehat{B}(0|\widehat{M}_{r})/v!),\mu\geq 1.$$ Here,
$q_{\mu,r}^{\theta_0}(x)$ is a polynomial which is determined by
the following conditions: when $x\in \Omega, $ $ q_{ \mu
,r}^{\theta_0}(x)=\sum\limits_{ j_1+\cdots+j_\kappa =\mu
}\prod\limits_{i=1}\limits^\kappa
\frac{s_{1+j_i}({\theta_0}^{-m_i})}{(j_i+1)!}\frac{1}{r^\kappa}D_{m_1}^{j_1}\cdots
D_{m_\kappa}^{j_\kappa}T(x|M_{{\theta_0}}), $ where
$s_0(x)=\frac{x-x^r}{x-1},s_j(x)=xs_{j-1}'(x),j\in {\bf Z}_+.$
\end{theorem}

In particular, when $M$ is a 1-prime matrix, a simple formulation
for $t(\cdot|M)$ is shown in the following theorem.

\begin{theorem}\cite{xu}\label{Th:form}
Under the condition of Theorem \ref{Theorem:DaMi,truncted power's
structure},when $M$ is a $1-$prime matrix,
$$
f_{\Omega}(\alpha |M)=p_{e,\Omega}(\alpha |M)+\sum_{{\theta} \in
A(M)\setminus e}{\theta} ^{\alpha
}\frac{1}{|det(M_{{\theta}})|}\prod_{w\in M\setminus
M_{{\theta}}}\frac{1}{1-{\theta}^{-w}}
1_{cone(M_{{\theta}})}(\Omega),
$$
where $p_{e,\Omega}(\alpha |M)$ is given in Theorem
\ref{Th:FormP}.
\end{theorem}

For the convenience of description, throughout the rest of the
paper, we suppose $M$ is a 1-prime matrix without further declaration.
 According to Theorem
\ref{Th:form}, to give a simple explicit formulation for $t({\bf
b}|M),$ we have to present  an explicit formulation for $T({\bf
x}|M)$. Moreover, to calculate the elements in $A(M)$ is a
non-trivial problem, hence, we have to predigest the
non-polynomial part in $t({\bf b}|M)$.

\section*{4. The generalized Popoviciu's formulation}
In this section, we are interested in $t({\bf n}|M)$, where
  $M=\begin{pmatrix}
  x_{1} & x_{2} & x_{3} \\
  y_{1} & y_{2} & y_{3}
\end{pmatrix}\in {\bf Z}^{2\times 3},$ ${\bf n}=(n_1,n_2)^T\in {\bf Z}_+^2.$
Without loss of generality, we suppose $
\frac{y_1}{x_1}<\frac{y_2}{x_2}< \frac{y_3}{x_3}$.  Obviously, for
the  matrix $M$, there exit two fundamental $M-$ cones, i.e.
$\Omega_1=\{(x,y)^T|(x,y)^T\in cone(M),\frac{y_1}{x_1}<
\frac{y}{x}< \frac{y_2}{x_2}\}$ and $\Omega_2=\{(x,y)^T|(x,y)^T\in
cone(M),\frac{y_2}{x_2}<\frac{y}{x}< \frac{y_3}{x_3}\}$ (See
Fig.1).

\begin{center}
\end{center}
\begin{center}
Fig.1.The fundamental M-cones.
\end{center}

To describe conveniently, we let $M_{ij}=\begin{pmatrix}
  x_{i} & x_{j}  \\
  y_{i} & y_{j}
\end{pmatrix}$, and let $Y_{ij}=det(M_{ij})$, where $i<j.$
To describe the explicit formulation for $t({\bf n}|M)$, we need
to define the fractional part function $\{x\}$ which denotes the
fractional part of $x$, i.e. $\{x\}=x-\lfloor x\rfloor$.

 In this section, our goal is to generalize the following beautiful formula due
to Popoviciu:
\begin{theorem}\label{th:popoviciu}\cite{popoviciu}
If $a$ and $b$ are relatively prime,
$$
t(n|(a,b))=\frac{n}{ab}-\{\frac{b^{-1}n}{a}\}-\{\frac{a^{-1}n}{b}\}+1,
$$
where $b^{-1}b\equiv 1\,\,\, mod \,\,\, a$, and $a^{-1}a\equiv
1\,\,\, mod \,\,\, b$, $n\in {\bf Z}_+$.
\end{theorem}

In order to generalize Theorem \ref{th:popoviciu}, we firstly
consider the explicit formulation for $T({\bf x}|M)$.

\begin{lemma}
Suppose the matrix
  $M=\begin{pmatrix}
  x_{1} & x_{2} & x_{3} \\
  y_{1} & y_{2} & y_{3}
\end{pmatrix}\in {\bf Z}^{2\times 3}.$
When ${\bf x}=(x,y)^T\in \overline{\Omega}_1, T({\bf x}|M
)=\frac{yx_1-xy_1}{(x_1y_2-y_1x_2)(x_1y_3-y_1x_3)}$; when ${\bf
x}=(x,y)^T\in \overline{\Omega}_2, T({\bf x}|M
)=\frac{xy_3-yx_3}{(x_2y_3-y_2x_3)(x_1y_3-y_1x_3)}.$
\end{lemma}
{\it proof:} Based on Theorem \ref{Th:numT} and $T({\bf
x}|M_{ij})=\frac{1}{det(M_{ij})}, {\bf x}\in cone(M_{ij}),i<j$,
the Lemma can be proved easily after a brief calculation. $\Box$

Hence, we obtain the conclusion as follows.

\begin{theorem}\label{Th:thfor}
Suppose the 1-prime matrix $M=\begin{pmatrix}
  x_{1} & x_{2} & x_{3} \\
  y_{1} & y_{2} & y_{3}
\end{pmatrix}$. When ${\bf n}=(n_1,n_2)^T\in
\overline{\Omega}_1\cap {\bf Z}^2$,
\begin{eqnarray*}
t({\bf n}|M)&=&\frac{n_2x_1-n_1y_1}{Y_{12}Y_{13}}-\{\frac{
(f_{12}Y_{13}+g_{12}Y_{23})^{-1}(n_2(f_{12}x_1+g_{12}x_2)-n_1(f_{12}y_1+g_{12}y_2))}{Y_{12}}\}\\
&-&\{\frac{
(f_{13}Y_{12}+g_{13}Y_{23})^{-1}(n_2(f_{13}x_1+g_{13}x_3)-n_1(f_{13}y_1+g_{13}y_3))}{Y_{13}}\}+1;
\end{eqnarray*}

when ${\bf n}=(n_1,n_2)^T\in \overline{\Omega}_2\cap {\bf Z}^2$,
\begin{eqnarray*}
t({\bf n}|M)&=&\frac{n_1y_3-n_2y_3}{Y_{23}Y_{13}}-\{\frac{
(f_{23}Y_{13}+g_{23}Y_{12})^{-1}(n_1(f_{23}x_3+g_{23}x_2)-n_2(f_{23}y_3+g_{23}y_2))}{Y_{23}}\}\\
&-&\{\frac{
(f_{13}Y_{12}+g_{13}Y_{23})^{-1}(n_1(f_{13}x_1+g_{13}x_3)-n_2(f_{13}y_1+g_{13}y_3))}{Y_{13}}\}+1,
\end{eqnarray*}
where, $f_{12},g_{12},f_{13}$,  $g_{13}$, $f_{23}$ and $g_{23}\in
{\bf Z}$ satisfy $gcd(f_{12}Y_{13}+g_{12}Y_{23},Y_{12})=1$
$gcd(f_{13}Y_{12}+g_{13}Y_{23},Y_{13})=1$ and
$gcd(f_{23}Y_{13}+g_{23}Y_{12},Y_{23})=1,$ moreover,
$(f_{12}Y_{13}+g_{12}Y_{23})^{-1}(f_{12}Y_{13}+g_{12}Y_{23})\equiv
1\,\,\, mod\,\,\, Y_{12},
(f_{13}Y_{12}+g_{13}Y_{23})^{-1}(f_{13}Y_{12}+g_{13}Y_{23})\equiv
1\,\,\, mod\,\,\, Y_{13},
(f_{23}Y_{13}+g_{23}Y_{12})^{-1}(f_{23}Y_{13}+g_{23}Y_{12})\equiv
1\,\,\, mod\,\,\, Y_{23}. $
\end{theorem}
{\it proof:} We only prove the case where $(n_1,n_2)^T\in
\overline{\Omega}_1\cap {\bf Z}^2.$ Based on Theorem
\ref{Th:FormP}, $p_{e,\Omega_1}({\bf x })$, which is the
polynomial part of $t(\cdot |M)$ on $\Omega_1$, is in the form of
$p_{0,\Omega_1}({\bf x})+p_{1,\Omega_1}({\bf x}),$ where for ${\bf
x}\in \Omega_1$, $p_{0,\Omega_1}({\bf x})=T({\bf x}|M),
p_{1,\Omega_1}({\bf x})=-(\sum_{|v|=1}D^vp_{0,\Omega_1}({\bf
x})(-i)D^v\widehat{B}(0 |M ))$. By the explicit formulation for
$T({\bf x}|M),$ we have $p_{0,\Omega_1}({\bf
x})=\frac{y_1x-x_1y}{(x_2y_1-y_2x_1)(x_1y_3-y_1x_3)}.$ After a
brief calculation, we have $p_{1,\Omega}({\bf
x})=\frac{1}{2}(\frac{1}{Y_{13}}+\frac{1}{Y_{12}}).$ Hence, the
polynomial part of $t({\bf n} |M)$ is
$\frac{n_2x_1-n_1y_1}{Y_{12}Y_{13}}+
\frac{1}{2}(\frac{1}{Y_{13}}+\frac{1}{Y_{12}})$. According to
Theorem \ref{Th:form}, we only need to consider  the sum
\begin{eqnarray*}
& &\sum_{\theta\in A(M)\setminus e}\frac{1}{|det(M_{\theta})|}
\prod_{w\in M\setminus M_{\theta}
}\frac{\theta^n}{1-\theta^{-w}}1_{cone(M_{\theta})}(\Omega_1)\\
& =& \frac{1}{Y_{12}}\sum_{\theta\in A(M)\setminus e\atop
M_\theta=M_{12}}\frac{\theta^n}{1-\theta^{-(x_3,y_3)}} +
\frac{1}{Y_{13}}\sum_{ { \theta \in A(M)\setminus e}\atop
{M_\theta=M_{13}} }\frac{\theta^n}{1-\theta^{-(x_2,y_2)}}.
\end{eqnarray*}

Recall $e^{\frac{2\pi i}{k}}$ is denoted by $W_k$. As pointed out
in \cite{DaMi on the solution}, the elements in the set
$\{\theta|\theta\in A(M),M_{\theta}=M_{12}\}$ have the form
$(W_{Y_{12}}^{\alpha_1^j},W_{Y_{12}}^{\alpha_2^j}),$ where
$(\alpha_1^j,\alpha_2^j)\in {\bf Z}^2,1\leq j\leq Y_{12}$.

 Consider firstly
\begin{eqnarray}\label{eq:formsum}
\frac{1}{Y_{12}}\sum_{\theta\in A(M)\setminus e
\atop M_\theta=M_{12}}\frac{\theta^{\bf n}}{1-\theta^{-(x_3,y_3)}}=
\frac{1}{Y_{12}}\sum_{j=1}^{Y_{12}-1}
\frac{W_{Y_{12}}^{n_1\alpha_1^j+n_2\alpha_2^j}}{1-W_{Y_{12}}^{-(x_3\alpha_1^j+y_3\alpha_2^j)}}.
\end{eqnarray}

We set $x_3\alpha_1^j+y_3\alpha_2^j\equiv k$  $mod$ $Y_{12}$.
Since $M$ is a 1-prime matrix,
$x_3\alpha_1^j+y_3\alpha_2^j\not\equiv
x_3\alpha_1^m+y_3\alpha_2^m$ $mod$ $Y_{12}$ when $j\neq m$. Hence,
$k$ runs over $[1,Y_{12}-1]\cap {\bf Z}.$

For $\theta\in \{\theta| \theta\in A(M),M_{\theta}=M_{12}\},$ we
have $\theta^{(x_1,y_1)}=\theta^{(x_2,y_2)}=1.$ Hence,
\begin{eqnarray}
x_1\alpha_1^j+y_1\alpha_2^j&\equiv &0 \,\,\, mod\,\,\, Y_{12}\label{eq:1},\\
x_2\alpha_1^j+y_2\alpha_2^j&\equiv &0 \,\,\, mod\,\,\, Y_{12}\label{eq:2},\\
x_3\alpha_1^j+y_3\alpha_2^j&\equiv & k \,\,\, mod\,\,\,
Y_{12}\label{eq:3}.
\end{eqnarray}

By $x_1$ on both sides of $(\ref{eq:3}),$ we have

\begin{eqnarray}\label{eq:eq1}
x_1x_3\alpha_1^j+x_1y_3\alpha_2^j\equiv x_1k \,\,\, mod\,\,\,
Y_{12}.
\end{eqnarray}
According to (\ref{eq:1}), we obtain
\begin{eqnarray}\label{eq:eq2}
x_1\alpha_1^j\equiv -y_1\alpha_2^j\,\,\, mod\,\,\, Y_{12}.
\end{eqnarray}

Substituting (\ref{eq:eq2}) into (\ref{eq:eq1}), we have
\begin{equation}\label{eq:diop1}
\alpha_2^jY_{13}\equiv x_1k \,\,\, mod\,\,\, Y_{12}.
\end{equation}
By using similar method,
\begin{equation}\label{eq:diop2}
\alpha_2^jY_{23}\equiv x_2k \,\,\, mod\,\,\, Y_{12}.
\end{equation}
Since $M$ is  a 1-prime matrix, there exits $f_{12},g_{12}\in {\bf
Z}$ such that  $gcd(f_{12}Y_{13}+g_{12}Y_{23},Y_{12})=1.$
Combining (\ref{eq:diop1}) and (\ref{eq:diop2}), we have
$$
\alpha_2^j(f_{12}Y_{13}+g_{12}Y_{23})\equiv
(f_{12}x_1+g_{12}x_2)k\,\,\, mod\,\,\, Y_{12}.
$$
Hence, $\alpha_2^j\equiv
(f_{12}Y_{13}+g_{12}Y_{23})^{-1}(f_{12}x_1+g_{12}x_2)k\,\,\, mod \,\,\, Y_{12}.$

Similarly,  $\alpha_1^j\equiv
-(f_{12}Y_{13}+g_{12}Y_{23})^{-1}(f_{12}y_1+g_{12}y_2)k\,\,\, mod
\,\,\, Y_{12}.$ Hence,(\ref{eq:formsum}) is reduced to
\begin{eqnarray}\label{eq:sums}
\frac{1}{Y_{12}}\sum_{k=1}^{Y_{12}-1}
\frac{W_{Y_{12}}^{(n_2(f_{12}x_1+g_{12}x_2)-n_1(f_{12}y_1+g_{12}y_2))(f_{12}Y_{13}+g_{12}Y_{23})^{-1}k}}
{1-W_{Y_{12
}}^{-k}}.
\end{eqnarray}
According to  discrete Fourier transforms,
\begin{equation}\label{eq:dis}
-\{\frac{t}{a}\}=\frac{1-a}{2a}+\frac{1}{a}\sum_{k=1}^{a-1}\frac{W_a^k}{1-W_a^{-k}},
\end{equation}

(\ref{eq:sums}) can be reduced to
\begin{eqnarray*}
\{\frac{(n_2(f_{12}x_1+g_{12}x_2)-n_1(f_{12}y_1+g_{12}y_2))(f_{12}Y_{13}
+g_{12}Y_{23})^{-1}}{Y_{12}}\} +\frac{1}{2}-\frac{1}{2Y_{12}}.
\end{eqnarray*}

Hence
\begin{eqnarray*}
 & &\frac{1}{Y_{12}}\sum_{\theta\in
A(M)\atop M_\theta=Y_{12}}\frac{\theta^n}{1-\theta^{-(x_3,y_3)}}\\
&=& -
\{\frac{(n_2(f_{12}x_1+g_{12}x_2)-n_1(f_{12}y_1+g_{12}y_2))(f_{12}Y_{13}+g_{12}Y_{23})^{-1}}{Y_{12}}\}
+\frac{1}{2}-\frac{1}{2Y_{12}}.
\end{eqnarray*}
By using similar method, we have
\begin{eqnarray*}
& &\frac{1}{Y_{13}}\sum_{\theta\in
A(M)\setminus e\atop M_\theta=Y_{13}}\frac{\theta^{\bf n}}{1-\theta^{-(x_2,y_2)}}\\
&=&-\{\frac{(n_2(f_{13}x_1+g_{13}x_3)-n_1(f_{13}y_1+g_{13}y_3))(f_{13}Y_{12}+g_{13}Y_{23})^{-1}}
{Y_{13}}\} +\frac{1}{2}-\frac{1}{2Y_{13}}.
\end{eqnarray*}

Hence, when $(n_1,n_2)^T\in v(\Omega_1|M)\cap {\bf Z}^2$,
\begin{eqnarray*}
t({\bf n}|M)&=&\frac{n_2x_1-n_1y_1}{Y_{12}Y_{13}}-\{\frac{
(f_{12}Y_{13}+g_{12}Y_{23})^{-1}(n_2(f_{12}x_1+g_{12}x_2)-n_1(f_{12}y_1+g_{12}y_2))}{Y_{12}}\}\\
&-&\{\frac{
(f_{13}Y_{12}+g_{13}Y_{23})^{-1}(n_2(f_{13}x_1+g_{13}x_3)-n_1(f_{13}y_1+g_{13}y_3))}{Y_{13}}\}+1.
\end{eqnarray*}
Note that $\overline{\Omega}_1\subset v(\Omega_1|M).$ Hence, when
${\bf n}\in \overline{\Omega}_1\cap {\bf Z}^2,$ the theorem holds.
 $\Box$

\begin{remark}
 If $f_{12}, g_{12},f_{13},g_{13},f_{23}$ and $g_{23}$ satisfy
$f_{12}Y_{23}+g_{12}Y_{13}=gcd(Y_{23},Y_{13}),f_{13}Y_{12}+g_{13}Y_{23}=gcd(Y_{12},Y_{23}),
f_{23}Y_{13}+g_{23}Y_{12}=gcd(Y_{13},Y_{12})$, then
$gcd(f_{12}Y_{13}+g_{12}Y_{23},Y_{12})=1,$
$gcd(f_{13}Y_{12}+g_{13}Y_{23},Y_{13})=1$ and
$gcd(f_{23}Y_{13}+g_{23}Y_{12},Y_{23})=1.$ Hence,one can determine
$f_{12}, g_{12},f_{13},g_{13},f_{23}$ and $g_{23}$ by Euclidean
algorithm. But in some special cases, such as $Y_{12},Y_{13}$ and
$Y_{23}$ are pairwise relative prime, there exits the simpler
method for obtaining them.
\end{remark}
\begin{corollary}\label{co:for}
Suppose $Y_{12},Y_{13}$ and $Y_{23}$ are pairwise relative prime.
When ${\bf n}=(n_1,n_2)^T\in \overline{\Omega}_1\cap {\bf Z}^2$,
\begin{eqnarray*}
t({\bf n}|M)=\frac{n_2x_1-n_1y_1}{Y_{12}Y_{13}}-\{\frac{
Y_{13}^{-1}(n_2x_1-n_1y_1)}{Y_{12}}\}-\{\frac{
Y_{12}^{-1}(n_2x_1-n_1y_1)}{Y_{13}}\}+1,
\end{eqnarray*}
where $Y_{13}^{-1}Y_{13}\equiv 1\,\,\, mod\,\,\, Y_{12}$ and
$Y_{12}^{-1}Y_{12}\equiv 1\,\,\, mod\,\,\, Y_{13}$. When ${\bf
n}=(n_1,n_2)^T\in \overline{\Omega}_2\cap {\bf Z}^2$,
\begin{eqnarray*}
t({\bf n}|M)=\frac{n_1y_3-n_2x_3}{Y_{23}Y_{13}}-\{\frac{
Y_{13}^{-1}(n_1x_3-n_2y_3)}{Y_{23}}\}-\{\frac{
Y_{23}^{-1}(n_1x_3-n_2y_3)}{Y_{13}}\}+1,
\end{eqnarray*}
where $Y_{13}^{-1}Y_{13}\equiv 1\,\,\, mod\,\,\, Y_{23}$ and
$Y_{23}^{-1}Y_{23}\equiv 1\,\,\, mod\,\,\, Y_{13}$.
\end{corollary}
{\it proof:} We firstly consider the case where ${\bf n}\in
\overline{\Omega}_1\cap {\bf Z}^2.$ Since $gcd(Y_{12},Y_{13})=1$,
$M$ is a 1-prime matrix. In Theorem \ref{Th:thfor}, we may set
$f_{12}=1,g_{12}=0,f_{13}=1,$ and $g_{13}=0.$ Hence, When ${\bf
n}=(n_1,n_2)^T\in \overline{\Omega}_1\cap {\bf Z}^2$,
\begin{eqnarray*}
t({\bf n}|M)=\frac{n_2x_1-n_1y_1}{Y_{12}Y_{13}}-\{\frac{
Y_{13}^{-1}(n_2x_1-n_1y_1)}{Y_{12}}\}-\{\frac{
Y_{12}^{-1}(n_2x_1-n_1y_1)}{Y_{13}}\}+1.
\end{eqnarray*}

Using similar method, when ${\bf n}=(n_1,n_2)^T\in
\overline{\Omega}_2\cap {\bf Z}^2$,
\begin{eqnarray*}
t({\bf n}|M)=\frac{n_1y_3-n_2x_3}{Y_{23}Y_{13}}-\{\frac{
Y_{13}^{-1}(n_1x_3-n_2y_3)}{Y_{23}}\}-\{\frac{
Y_{23}^{-1}(n_1x_3-n_2y_3)}{Y_{13}}\}+1.
\end{eqnarray*}
 $\Box$

\begin{remark}
An interesting observation is that the formulation presented in
Corollary \ref{co:co1} is remarkably similar with Popoviciu's
formulation.
\end{remark}

We now turn to consider the special  case where
$\frac{y_1}{x_1}=\frac{y_2}{x_2}$. Without loss of generality, we
suppose $M=\begin{pmatrix}
  kx_{1} & lx_{1} & x_{3} \\
  ky_{1} & ly_{1} & y_{3}
\end{pmatrix}$.
In this case, there exits only one fundamental M-cone, which is
denoted as $\Omega$. Moreover, since $M$ is a 1-prime matrix, we
have $gcd(k,l)=1,x_1y_3-y_1x_3=1.$
 Then we have
\begin{theorem}\label{th:spM}
Suppose $\frac{y_1}{x_1}<\frac{y_3}{x_3}.$ When $M=\begin{pmatrix}
  kx_{1} & lx_{1} & x_{3} \\
  ky_{1} & ly_{1} & y_{3}
\end{pmatrix},$
$t({\bf n
}|M)=\frac{x_3n_2-y_3n_1}{kl}-\{\frac{l^{-1}}{k}(n_1y_3-n_2x_3)\}-
\{\frac{k^{-1}}{l}(n_1y_3-n_2x_3)\}+1,$ where ${\bf
n}=(n_1,n_2)^T\in \overline{\Omega}\cap {\bf Z}^2.$
\end{theorem}
{\it proof:} By using the recurrence formulation for $T({\bf x}|M)$, we have
 $T({\bf x}|M)=\frac{x_3y-y_3x}{kl}.$ Hence, the polynomial
part of $t(\cdot|M)$ is $\frac{x_3y-y_3x}{kl}+\frac{1}{2}(\frac{1}{k}+
\frac{1}{l}).$ We
now only need to consider the sums
$$
\frac{1}{k}\sum_{\theta:M_{\theta}=Y_k}\frac{\theta^n}{1-\theta^{-(lx_1,ly_1)}},
\frac{1}{l}\sum_{\theta:M_{\theta}=Y_l}\frac{\theta^n}{1-\theta^{-(kx_1,ky_1)}}.
$$
By using the similar method with the one presented in the proof of
Theorem \ref{Th:thfor}, we have $
\frac{1}{k}\sum\limits_{\theta\in A(M) \setminus e \atop
M_{\theta}=Y_k}\frac{\theta^n}{1-\theta^{-(lx_1,ly_1)}}
=-\{l^{-1}\frac{n_1y_3-n_2x_3}{k}\}+\frac{1}{2}-\frac{1}{2k},$
$\frac{1}{l}\sum\limits_{\theta\in A(M)\atop
M_{\theta}=Y_l}\frac{\theta^n}{1-\theta^{-(kx_1,ky_1)}}
=-\{k^{-1}\frac{n_1y_3-n_2x_3}{k}\}+\frac{1}{2}-\frac{1}{2l}.$
Note that $\overline{\Omega}\subset v(\Omega|M)$. By Theorem
\ref{Th:form}, when ${\bf n}=(n_1,n_2)^T\in \overline{\Omega},$
\begin{eqnarray*}
& &t({\bf n }|M)=\frac{x_3y-y_3x}{kl}+\frac{1}{2}(\frac{1}{k}+
\frac{1}{l})+
\frac{1}{k}\sum\limits_{\theta:M_{\theta}=Y_k}\frac{\theta^n}{1-\theta^{-(lx_1,ly_1)}}
+\frac{1}{l}\sum\limits_{\theta:M_{\theta}=Y_l}\frac{\theta^n}{1-\theta^{-(kx_1,ky_1)}}\\
&=&\frac{x_3n_2-y_3n_1}{kl}-\{\frac{l^{-1}}{k}(n_1y_3-n_2x_3)\}-
\{\frac{k^{-1}}{l}(n_1y_3-n_2x_3)\}+1.
\end{eqnarray*}

$\Box$

\begin{remark}
When the matrix $M$ is of the form $\begin{pmatrix}
  x_{1} & kx_{2} & lx_{2} \\
  y_{1} & ky_{2} & ly_{2}
\end{pmatrix},$ a similar result can be obtained using the same method
with the one presented in Theorem \ref{th:spM}.
\end{remark}

\section*{5.Linear Diophantine problem of Frobenius}

Consider the linear Diophantine equation
\begin{equation}\label{eq:dio}
x_1a_1+\cdots x_na_n=N,
\end{equation}
where, $a_i\in {\bf Z}_+, gcd(a_1,\cdots,a_n)=1.$

It is well known that for all sufficiently large $N$ the equation
has solutions. The Frobenius problems asks us to find the largest
integer for which no solution exists. We call the largest integer
the Frobenius number and denote it by $f(a_1,\cdots,a_n).$ For
$n=2$ the largest $N$ for which no solution exists can be
explicitly written as $a_1a_2-a_1-a_2,$ i.e.
$f(a_1,a_2)=a_1a_2-a_1-a_2.$ But no such formula exists for $n\geq
3.$

As pointed out in \cite{xu}, when $gcd\{|Y|:Y\in {\mathcal
B}(M)\}=1,$ for all sufficiently large $N$ the linear Diophantine
equations $Mx=N{\bf n}$ has solution, where ${\bf n}\in cone(M).$
Naturally, we hope to find the largest integer $N$ for which no
solution exits, which is denoted as $f(M,{\bf n}).$ In particular,
we are interested in the linear Diophantine equations $M_0x=N{\bf
n},$ where $ M_0=\begin{pmatrix}
  x_{1} & x_{2} & x_{3} \\
  y_{1} & y_{2} & y_{3}
\end{pmatrix}, {\bf n}\in cone(M_0).
$ In fact, the generalized Frobenius number $f(M_0,{\bf n})$ is a
generalization of $f(a_1,a_2).$

Recall $M_{ij}=\begin{pmatrix}
  x_{i} & x_{j}  \\
  y_{i} & y_{j}
\end{pmatrix}$ and  $Y_{ij}=det(M_{ij}).$
In the following theorem, we shall present an upper  boundary for
$f(M_0,{\bf n}).$
\begin{theorem}\label{th:genf}
 Suppose $Y_{12},Y_{13}$ and $Y_{23}$ are pairwise relative prime. For ${\bf n}\in \overline{\Omega}_1\cap {\bf Z}^2,$
 $f(M_0,{\bf n})< \frac{Y_{12}Y_{13}-Y_{12}-Y_{13}+1}{n_2x_1-n_1y_1}.$
 For ${\bf n}\in \overline{\Omega}_2\cap {\bf Z}^2,$
 $f(M_0,{\bf n})< \frac{Y_{23}Y_{13}-Y_{23}-Y_{13}+1}{n_1y_3-n_2x_3}.$
\end{theorem}
{\it proof:} We only prove the case where ${\bf n}\in
\overline{\Omega}_1\cap {\bf Z}^2.$ Note $t(N{\bf
n}|M)=\frac{N(n_2x_1-n_1y_1)}{Y_{12}Y_{13}}-\{\frac{
(Y_{13})^{-1}(N(n_2x_1-n_1y_1))}{Y_{12}}\}-\{\frac{
(Y_{12})^{-1}(N(n_2x_1-n_1y_1))}{Y_{13}}\}+1=t(N(n_2x_1-n_1y_1)|(Y_{12},Y_{13})).$
Since when $N(n_2x_1-n_1y_1)\geq Y_{12}Y_{13}-Y_{12}-Y_{13}+1,$ $
t(N(n_2x_1-n_1y_1)|Y_{12},Y_{13})=t(N{\bf n}|M)>0.$ Hence, when
$N\geq \frac{Y_{12}Y_{13}-Y_{12}-Y_{13}+1}{(n_2x_1-n_1y_1)},$
$t(N{\bf n}|M)>0.$ So, $f(M,{\bf n})<
\frac{Y_{12}Y_{13}-Y_{12}-Y_{13}+1}{n_2x_1-n_1y_1}.$ $\Box$

\begin{remark}
Theorem \ref{th:genf} only gives an upper boundary for $f(M_0,{\bf
n}).$ According to the proof of Theorem \ref{th:genf}, giving the
exact value of $f(M_0,{\bf n})$ is equivalent for any given
$b_0\in {\bf Z }$ determining the largest integer $N$ for which
the Diophantine equation $x_1a_1+x_2a_2=Nb_0$ no solution exits.
\end{remark}

\section*{6 Two-dimension vector partition functions}

We now turn to the general case. Let $M=\begin{pmatrix}
  x_{1} & x_{2} & \cdots &x_n \\
  y_{1} & y_{2} & \cdots & y_n
\end{pmatrix}$ be  a $2\times n$ integer matrix.
and $\frac{y_{i-1}}{x_{i-1}}<\frac{y_i}{x_i},i=2,\cdots,n$.

For the matrix $M$, there exist $n-1$ fundamental $M-$ cones.
Denote them as $\Omega_i:=\{(x,y)^T|(x,y)^T\in
cone(M),\frac{y_i}{x_i}<\frac{y}{x}<\frac{y_{i+1}}{x_{i+1}}\},i=1,\cdots,n-1$
respectively. In this section, we shall discuss the explicit
formulation for $t({\bf b}|M)$. First, we present an explicit
formulation for $T({\bf x}|M)$.
\begin{theorem}\label{Th:2dT} For ${\bf x}=(x,y)^T\in {\bf R}^2,$
$$
T({\bf x
}|M)=\frac{1}{(n-2)!}\sum_{i=1}^n\frac{(y_ix-x_iy)_+^{n-2}}{\prod_{j\neq
i}(y_ix_j-y_jx_i) },
$$
where, $(y_ix-x_iy)_+=\begin{cases}
y_ix-x_iy , & y_ix-x_iy\geq 0,\\
0 , & otherwise.
\end{cases}
$
\end{theorem}
{\it proof:} According to the definition of $(y_ix-x_iy)_+$ we
only need to prove that when ${\bf x}\in \Omega_k$,
$T({\bf x}|M)=\frac{1}{(n-2)!}\sum_{i=k+1}^n\frac{(y_ix-x_iy)^{n-2}}{\prod_{j\neq
i}(y_ix_j-y_jx_i) }.$

 We argue by induction on $n$.
Initially, when $n=2,3$ the theorem certainly holds. In the
inductive step,  we assume that when $n=n_0$ the theorem holds and
we consider the case when $n=n_0+1$.

 According to the definition of $(y_ix-x_iy)_+$ we
only need to prove that for ${\bf x}\in \Omega_k,T({\bf x
}|M)=\frac{1}{(n_0-1)!}\sum\limits_{i=k+1}^{n_0+1}\frac{(y_ix-x_iy)^{n_0-1}}{\prod_{j\neq
i}(y_ix_j-y_jx_i) },$
 where $M$ is a $2\times (n_0+1)$ matrix.

 After a brief calculation,  it is easy for
obtaining $ {\bf x
}=\frac{xy_{k+1}-x_{k+1}y}{y_{k+1}x_k-y_kx_{k+1}}(x_k,y_k)^T+
\frac{xy_{k}-x_{k}y}{y_{k+1}x_k-y_kx_{k+1}}(x_{k+1},y_{k+1})^T.$
Based on the recurrence formulation of $T(\cdot |M)$ , we have
$$T({\bf x}|M)=
\frac{1}{n_0-1}(\frac{xy_{k+1}-x_{k+1}y}{y_{k+1}x_k-y_kx_{k+1}}
T({\bf x}|M\setminus (x_k,y_k)^T)+
\frac{xy_{k}-x_{k}y}{y_{k+1}x_k-y_kx_{k+1}} T({\bf x}|M\setminus
(x_{k+1},y_{k+1})^T)).$$ By the inductive hypothesis, $ T({\bf
x}|M\setminus
(x_k,y_k)^T)=\frac{1}{(n_0-2)!}\sum\limits_{i=k+1}^{n_0+1}\frac{(y_ix-x_iy)^{n_0-2}}{\prod\limits_{j\neq
i,j\neq k}(y_ix_j-y_jx_i) },$ $T({\bf x}|M\setminus
(x_{k+1},y_{k+1})^T)=\frac{1}{(n_0-2)!}\sum\limits_{i=k+2}^{n_0+1}\frac{(y_ix-x_iy)^{n_0-2}}{\prod\limits_{j\neq
i}(y_ix_j-y_jx_i) }.$ Then we obtain
\begin{eqnarray*}
& &T({\bf
x}|M)=\frac{1}{(n_0-1)!}(\frac{xy_{k+1}-x_{k+1}y}{y_{k+1}x_k-y_kx_{k+1}}\sum_{i=k+1}^{n_0+1}
\frac{(y_ix-x_iy)^{n_0-2}(x_ky_i-y_kx_i)}{\prod\limits_{j\neq
i}(y_ix_j-y_jx_i)} \\
&+&\frac{xy_{k}-x_{k}y}{y_{k+1}x_k-y_kx_{k+1}}\sum_{i=k+2}^{n_0+1}
\frac{(y_ix-x_iy)^{n_0-2}(x_{k+1}y_i-y_{k+1}x_i)}{\prod\limits_{j\neq
i}(y_ix_j-y_jx_i)} )\\
&=&\frac{1}{(n_0-1)!}(\frac{xy_{k+1}-yx_{k+1}}{y_{k+1}x_k-y_kx_{k+1}}\frac{(y_{k+1}x-x_{k+1}y)^{n_0-2}
(x_ky_{k+1}-y_kx_{k+1})}{\prod\limits_{j\neq i}(y_ix_j-y_jx_i)}
+\frac{1}{y_{k+1}x_k-y_kx_{k+1}}\\
& &\sum_{i=k+2}^{n_0+1}
(y_ix-x_iy)^{n_0-2}(\frac{(xy_{k+1}-x_{k+1}y)(x_ky_i-y_kx_i)-
(xy_{k}-x_{k}y)(x_{k+1}y_i-y_{k+1}x_i)}{\prod\limits_{j\neq
i}(y_ix_j-y_jx_i)}))\\
&=&\frac{1}{(n_0-1)!}\sum_{i=k+1}^{n_0+1}\frac{(y_ix-x_iy)^{n_0-1}}{\prod\limits_{j\neq
i}(y_ix_j-y_jx_i) }.
\end{eqnarray*}
Thus, when $n=n_0+1$ the theorem holds also, which completes the
inductive step and the proof. $\Box$

The following statements follow from Theorem \ref{Th:2dT}.

\begin{corollary}\label{co:co1}
$$
D^{v_1,v_2}T({\bf
x}|M)=\frac{1}{(n-2-v_1-v_2)!}\sum_{i=1}^n\frac{(y_ix-x_iy)_+^{n-2-v_1-v_2}}{\prod\limits_{j\neq
i }(y_ix_j-y_jx_i)}y_i^{v_1}(-x_i)^{v_2}.
$$
\end{corollary}

We now turn to non-polynomial part in $t(\cdot |M).$

In \cite{BeckFrobenius}, the Fourier-Dedekind sum is defined as
$\sigma_t(C;n)=\frac{1}{n}\sum\limits_{\lambda^n=1\neq\lambda}\frac{\lambda^t}{\prod_{c\in
C }(\lambda^c-1)}$, where $C$ is an integer multiset and $n$ is an
integer. To simplify the non-polynomial part in $t(\cdot |M),$ we
naturally arrived at the sums
\begin{equation}\label{eq:dedekind}
\frac{1}{Y_{ij}}\sum_{\theta^{M_{ij}}=1,\theta\neq e
}\theta^n\prod_{\omega\in M\setminus
M_{ij}}\frac{1}{1-\theta^{-\omega}},
\end{equation}
 which is considered as a generalized Fourier-Dedekind sum.
Here, $\theta ^{M_{ij}}=1$ means $\theta^{m}=1$ for any $m\in
M_{ij}.$
 In fact, it is a non-trivial problem
for computing all complex vectors satisfying $\theta^{M_{ij}}=1.$
In the following Lemma, we shows the generalized Fourier-Dedekind
sums (\ref{eq:dedekind}) can be converted into the 1-dimensional
Fourier-Dedekind sums.
\begin{lemma}\label{Le:Dedekind}
When $M$ is a 1-prime matrix, for any given integer $m,$ $1\leq
m\leq n, m\neq i,j,$
$$
\frac{1}{Y_{ij}}\sum_{\theta^{M_{ij}}=1\atop \theta\neq e
}\theta^n\prod_{\omega\in M\setminus
M_{ij}}\frac{1}{1-\theta^{-\omega}}
=\sigma_{t_{ij}}(C_{ij};Y_{ij}),
$$
where $C_{ij}=\cup_{1\leq h\leq n,h\neq i,h\neq j
}\{(fY_{im}+gY_{jm})^{-1}(-(fy_i+gy_j)x_h+(fx_i+gx_j)y_h)
\},t_{ij}=(fY_{im}+gY_{jm})^{-1}(-(fy_i+gy_j)n_1+(fx_i+gx_j)n_2)+\sum_{c\in
C_{ij}}c, $ where $f,g\in {\bf Z}$ satisfy
$gcd(fY_{im}+gY_{jm},Y_{ij})=1,$  moreover,
$(fY_{im}+gY_{jm})^{-1}(fY_{im}+gY_{jm})\equiv 1, mod \,\,\,
Y_{ij}.$
\end{lemma}
{\it proof:}
As pointed out in \cite{DaMi on the solution}, the elements in the
set $\{\theta|\theta\in A(M),M_{\theta}=M_{ij}\}$ have the form
$(W_{Y_{ij}}^{\alpha_1^l},W_{Y_{12}}^{\alpha_2^l}),$ where
$(\alpha_1^l,\alpha_2^l)\in {\bf Z}^2,1\leq l\leq Y_{ij}$.

Hence,
\begin{eqnarray}\label{eq:sum1}
\frac{1}{Y_{ij}}\sum_{\theta^{M_{ij}}=1\atop \theta\neq e
}\theta^n\prod_{\omega\in M\setminus
M_{ij}}\frac{1}{1-\theta^{-\omega}}=
\frac{1}{Y_{ij}}\sum_{l=1}^{Y_{ij}-1}
\frac{W_{Y_{ij}}^{n_1\alpha_1^l+n_2\alpha_2^l}}{\prod\limits_{h\neq
i,h\neq j} (1-W_{Y_{ij}}^{-(x_h\alpha_1^l+y_h\alpha_2^l)})}.
\end{eqnarray}

Noting $m\neq i,m\neq j,$ we set
$x_m\alpha_1^l+y_m\alpha_2^l\equiv k$ $mod$ $Y_{ij}$. Since $M$ is
a 1-prime matrix, $k$ runs over $[1,Y_{ij}-1]\cap {\bf Z}.$ Using
the similar method with the one in the proof of Theorem
\ref{Th:thfor}, we have
\begin{eqnarray*}
\alpha_1^l&\equiv &-(f_{ij}Y_{im}+g_{ij}Y_{jm})^{-1}(f_{ij}y_i+g_{ij}y_j)k \,\,\, mod\,\,\, Y_{ij},\\
\alpha_2^l&\equiv
&(f_{ij}Y_{im}+g_{ij}Y_{jm})^{-1}(f_{ij}x_i+g_{ij}x_j)k \,\,\,
mod\,\,\, Y_{ij}.
\end{eqnarray*}

Hence, (\ref{eq:sum1}) is converted into
\begin{eqnarray*}
& &\frac{1}{Y_{ij}}\sum_{k=1}^{Y_{ij}-1}
\frac{W_{Y_{ij}}^{(n_2(f_{ij}x_i+g_{ij}x_j)-n_1(f_{ij}y_i+g_{12}y_j))(f_{ij}Y_{im}+g_{ij}Y_{jm})^{-1}k}}
{\prod\limits_{h\neq i,h\neq j}(1-W_{Y_{ij
}}^{-(f_{ij}Y_{im}+g_{ij}Y_{jm})^{-1}(-x_h(f_{ij}y_i+g_{ij}y_j)+y_h(f_{ij}x_i+g_{ij}x_j))k})}\\
&=&\sigma_{t_{ij}}(C_{ij};Y_{ij}).
\end{eqnarray*}
$\Box$

\begin{remark}
When $|Y_{ij}|=1$,  since $\{\theta:\theta^{M_{ij}}\}=\{e\},$ the
terms in $\sigma_{t_{ij}}(C_{ij}:Y_{ij})$ disappear.
\end{remark}

Combining Theorem \ref{Th:FormP}, Theorem \ref{Th:form}, Theorem
\ref{Th:2dT} and Lemma \ref{Le:Dedekind}, we can present a
simplified  formulation for $t(\cdot |M)$.
\begin{theorem}\label{th:2dvfor}
Suppose $M=\begin{pmatrix}
  x_{1} & x_{2} & \cdots &x_n \\
  y_{1} & y_{2} & \cdots & y_n
\end{pmatrix}$ is  a $2\times n$ integer 1-prime matrix and
$\frac{y_i}{x_i}<\frac{y_{i+1}}{x_{i+1}}$. When ${\bf
n}=(n_1,n_2)^T\in {\Omega}_k\cap {\bf Z}^2$,
$$
 t({\bf n}|M)=p_{e,\Omega_k}({\bf n})+\sum_{(i,j)\in \{(i,j): i\leq k<j \} }\sigma_{t_{ij}}(C_{ij};Y_{ij}),
$$
where,
$p_{e,\Omega_k}({\bf x})=\sum_{j=0}^{n-2}p_{j,\Omega_k}({\bf x}),
p_{0,\Omega_k}({\bf x})=\frac{1}{n-2}
\sum_{l=k+1}^n\frac{(y_lx-x_ly)^{n-2}}{\prod_{j\neq
l}(y_lx_j-y_jx_l)},
p_{j,\Omega_k}({\bf x})=-\sum_{l=0}^{j-1}(\sum_{|v|=j-l}D^vp_{l,\Omega_k}({\bf x})
(-i)^{|v|}\frac{D^v\widehat{B}(0|M)}{v!}),$ $t_{ij}$ and $C_{ij}$  are defined in Lemma
\ref{Le:Dedekind}.
\end{theorem}
{\it proof:} Based on Theorem \ref{Th:form}, when ${\bf n}\in \Omega_k\cap {\bf Z}^2,$
$$
t({\bf n}|M)=p_{e,\Omega_k}({\bf n})+\sum_{{\theta} \in
A(M)\setminus e}{\theta} ^{\bf n
}\frac{1}{|det(M_{{\theta}})|}\prod_{w\in M\setminus
M_{{\theta}}}\frac{1}{1-{\theta}^{-w}}
1_{cone(M_{{\theta}})}(\Omega_k),
$$
where, the $p_{e,\Omega_k}$ can be determined easily. Since $M$ is a 1-prime matrix,
\begin{eqnarray*}
& &\sum_{{\theta} \in A(M)\setminus e}{\theta} ^{\bf n
}\frac{1}{|det(M_{{\theta}})|}\prod_{w\in M\setminus
M_{{\theta}}}\frac{1}{1-{\theta}^{-w}}
1_{cone(M_{{\theta}})}(\Omega_k)\\
&=&\sum_{i<j} \frac{1}{Y_{ij}}\sum_{\theta^{M_{ij}}=1\atop
\theta\neq e }\theta^{\bf n}\prod_{\omega\in M\setminus
M_{ij}}\frac{1}{1-\theta^{-\omega}}1_{cone(M_{ij})}(\Omega_k).
\end{eqnarray*}
Based on Lemma \ref{Le:Dedekind}, the above sum becomes as follows:
\begin{equation}\label{eq:sum}
\sum_{i<j}\sigma_{t_{ij}}(C_{ij}:Y_{ij})1_{cone(M_{ij})}(\Omega_k).
\end{equation}
Since when $k\geq j$ or $k< i$, $cone(M_{ij})\cap \Omega_k=\emptyset$. Hence, (\ref{eq:sum})
is converted into
\begin{equation}
\sum_{(i,j)\in \{(i,j): i\leq k<j \}
}\sigma_{t_{ij}}(C_{ij}:Y_{ij}).
\end{equation}
The theorem holds.
$\Box$

The explicit formulation presented in Theorem \ref{th:2dvfor}
contains $D^v\widehat{B}(0|M).$ Note $$ \widehat{B}(\zeta
|M)=\prod_{j=1}^n\frac{1-exp(-i\zeta^Tm_j)}{i\zeta^Tm_j},\zeta \in
{\bf C}^s.
$$

The following assertion is obvious.
$$
D^{v_1,v_2}\widehat{B}(0|M)=(-i)^{v_1+v_2}\sum_{k_1+\cdots+k_n=v_1}\sum_{l_1+\cdots+l_n=v_2}\frac{v_1!}
{k_1!\cdots k_n!}\frac{v_2!}{l_1!\cdots
l_n!}\prod_{j=1}^n\frac{x_j^{k_j}y_j^{l_j}}{k_j+l_j+1}.
$$

 Using  Theorem \ref{th:2dvfor}, we shall present an explicit formulation for an
actual vector partition function, which is the same with the one
presented in \cite{dgbeck}. By using Theorem \ref{th:2dvfor}, it
is indeed easier for obtaining the explicit formulation for the
actual vector partition function.

\begin{example}\label{ex:ex1}
Let  $A=\begin{pmatrix}
1 & 2 & 1 &0 \\
0 & 1 & 1 &1
\end{pmatrix}.$ We denote by $A_{ij}$ the square matrix containing
the $i$th and the $j$th columns in $A$.

For the matrix $A$, there exit three fundamental cones, which are
denoted  as $\Omega_1,\Omega_2$ and $\Omega_3$ respectively. We
shall discuss the explicit formulation for $t({\bf n}|A).$ After a
brief calculation, we have

$
 T({\bf x}|A)=
\begin{cases}
\frac{y^2}{2} ,& {\bf x}\in \Omega_1,\cr
\frac{1}{4}(-x^2+4xy-2y^2),&{\bf x}\in \Omega_2,\cr
\frac{x^2}{4},&{\bf x}\in \Omega_3.
\end{cases}
$

 Hence,
$p_{0,\Omega_1}=\frac{y^2}{2}$. According to Theorem \ref{Th:2dT},
$p_{1,\Omega_1}=3/2y$ and $p_{2,\Omega_1}=1$ respectively. Since
for any $1<j\leq 3,|det(Y_{1j})|=1$, the terms in Fourier-Dedekind
sum shall not appear when ${\bf n}\in \Omega_1\cap {\bf Z}^2$.
 Based on Theorem \ref{Th:2dT}, we have when
${\bf n}\in \Omega_1\cap {\bf Z}^2 ,t({\bf n|A
})=\frac{n_2^2}{2}+\frac{3n_2}{2}+1.$

Similarly, $p_{0,\Omega_2}=\frac{1}{4}(-x^2+4xy-2y^2),
p_{1,\Omega_2}=\frac{x+y}{2},p_{2,\Omega_2}=\frac{7}{8}$. Based on
Lemma \ref{Le:Dedekind}, the non-polynomial part is
$\frac{1}{Y_{23}}\sum\limits_{\theta^{A_{23}}=1,\theta\neq e
}\theta^n\prod\limits_{\omega\in A\setminus
A_{ij}}\frac{1}{1-\theta^{-\omega}}=(-1)^{n_1}$. Hence, when ${\bf
n}\in \Omega_2\cap {\bf Z}^2 , t({\bf n}|A)=
n_1n_2-\frac{n_1^2}{4}-\frac{n_2^2}{2}+\frac{n_1+n_2}{2}+\frac{7}{8}+\frac{(-1)^{n_1}}{8}.$

Using the same method with the above,  we obtain
$p_{0,\Omega_3}=\frac{x^2}{4},p_{1,\Omega_3}=x,p_{2,\Omega_3}=\frac{7}{8}.$

Hence,  $ t({\bf n}|A)=
\begin{cases}
\frac{n_2^2}{2}+\frac{3n_2}{2}+1,& {\bf n}\in \Omega_1\cap {\bf
Z}^2 \\
n_1n_2-\frac{n_1^2}{4}-\frac{n_2^2}{2}+\frac{n_1+n_2}{2}+\frac{7}{8}+\frac{(-1)^{n_1}}{8},&{\bf n}\in \Omega_2\cap {\bf Z}^2\\
\frac{n_1^2}{4}+n_1+\frac{7}{8}+\frac{(-1)^{n_1}}{8},&{\bf n}\in
\Omega_3\cap {\bf Z}^2.
\end{cases}
$
\end{example}

\begin{remark}
In Theorem  \ref{th:2dvfor}, when the case of
$\frac{y_i}{x_i}=\frac{y_j}{x_j}$ happens, the explicit
formulation for $T({\bf x}|M)$ can be obtained by taking the
limit. Using similar method with the one in the proof of Theorem
\ref{th:spM}, an explicit formulation for $t({\bf n}|M)$ can be
given also.
\end{remark}

\begin{remark}
 To simplify any-dimensional vector partition functions, we have to
give an explicit formulation for multivariate truncated power
functions $T({\bf x}|M)$ and compute the chamber complex
consisting of the fundamental $M-$cones, which are indeed
challenging problems.
\end{remark}

% \section{}
% \label{}


\begin{thebibliography}{999}
\small



\bibitem{BeckFrobenius}  M. Beck, R. Diaz and S. Robins, The Frobenius Problem,
Rational Polytopes, and Fourier-Dedekind Sums, J. Number
Theory.,(96) (2002), 1-21.

\bibitem{beckbook} M. Beck, and S. Robins, Computing the continuous
discretely integer point enumeration in polyhedron, to appear in
Springer Undergraduater Texts in Mathematics.

\bibitem{volume1} B. Bueler, A. Enge and K. Fukuda, Exact volume
computation for polytopes: A practical study. In:
Polytopes---Combinatorics and Computation, G. Kalai, and G. M.
Ziegler, Eds., Birh\"{a}user Verlag, Basel,(2000).



\bibitem{dgbeck} M. Beck, The partial-Fractiona method for counting
solutions to integral linear systems, Discrete Comput. Geom.
32:437-446(2004).

\bibitem{cappell} S.E. Caappell and J.L. Shaneson, Genera of
algebraic varieties and counting of lattice points, Bull.
A.M.S.30(1994),62-69.


\bibitem{dahemenpower}  W.Dahmen, On multivariate B-splines, SIAM J. Numer. Anal.
17(1980),179-191.

\bibitem{DahmenDM}  W.Dahmen and C.A. Micchelli, Translates of multivariate
splines, linear Algebra Appl. 52/53,(1983),217-234.

\bibitem{dahmen1}  W.Dahmen, and C.A. Micchelli, Recent progress in
multivariate splines, in {\it Approximation Theory IV}(C. K. Chui,
L. L. Schumaker, and J. Ward,Eds.), Academic, New
York,1983,pp.27-121.

\bibitem{DaMi on the solution}  W.Dahmen and C.A. Micchelli, On the solution of certain
systems of partial difference equations and linear dependence of
translates of box splines, Trans. Amer. Math. Soc.
292,(1985),305-320.


\bibitem{dahmen on the local linear}  W.Dahmen and C.A.Micchelli,On the local linear independence
of translate of a box spline, Studia Math. 82,(1985),243-263.

\bibitem{dahmen2}  W.Dahmen and C.A. Micchelli, The number of solutions to
linear diophantine equations and multivariate splines, Trans.
Amer. Math. Soc.308,(1988),509-532.


\bibitem{deboor1}  C.de Boor and R.Devore,
Approximation by smooth multivariate splines, Trans. Amer. Math.
Soc.276(1983),775-788.

\bibitem{deboor2}  C.de Boor and K.H\"{o}llig, B-splines from parallelepipeds,
J. Analyse Math,42(1982/83).

\bibitem{Astati} Dinwoodie, I.,Stochastic simulation on integer
constraint sets, SIAM J. Optimization 9(1999)53-61.

\bibitem{deboorbook} C.de Boor, K.H\"{o}llig  and S.Riemenschneider, { Box Splines},
 Springer-Verlag, New York,(1993).

\bibitem{deloera} Jes\'{u}s A. De Loera and Bernd Sturmfels,
Algebraic unimodular counting, Math. Program., Ser. B 96 (2003) 2,
183-203.


\bibitem{stanley}  R. Stanley, Enumerative Combinctorics, Vol.1, Wadsworth,
Belmont, Calif.,1986.

\bibitem{Jia:discret spline}  R.Q. Jia, Multivariate discrete splines and linear
diophantine equations, Trans. Amer. Math. Soc.340:1(1993),179-197.

\bibitem{jiamagic}  R.Q. Jia, Symmetric magic squares and multivariate
splines,Linear Algebra Appl.,250(1997)69-103.



\bibitem{numeircalT} C.A.Micchelli, On a numerically efficient
method for computing multivariate B-splines,in Multivariate
Approximation Theory, W. Schempp, and K.Zeller, eds.,
Birkh\"{a}user, Basel,(1979),211-248.

\bibitem{number} Nijehuis, A. and Wilf, H., Representation of
integers by linear forms in nonegative integers, J. Number Theory
4(1972),98-106.

\bibitem{ehrhart}  E. Ehrhart, Sur un probl\`{e}me de g\'{e}om\'{e}trie
diophantienne lin\'{e}aire II, J. Reine Angew. Math.
227(1967),25-49.

\bibitem{wangxu} R.H. Wang and Z.Q. Xu, Multivariate splines and lattice
points in rational polytope, J. Comp. Appl. Math.,159 (2003)
149-159.


\bibitem{pommersheim} J. Pommersheim, Toric varieties, lattice
points, and Dedekind sums, Math. Ann. 295(1993),1-24.

\bibitem{kantor} J.M. Kantor and A. Khovanskii, Une application du
Th\'{e}or\`{e}me de Riemann-Roch combinatoire au polyn\^{o}me
d'Ehrhar des polytopes entier de $R^d,$ C.R. Acad. Sci. Paris,
Series I 317 (1993),501-507.

\bibitem{gfrobenius} A. Rycerz, Conductors and the generalizaed
problem of Frobenius, Discussiones Mathematicae 14(1994)15-20.

\bibitem{rando} Welsh, D. Approximate counting , in Surveys in
Combinatorics, edited by R.A. Bailey, London Mathematical Society
Lecture Notes, Vol.241,1997.

\bibitem{repre} Schmidt, J.R. and Bincer, A. The Kostant partition
function for simple Lie algebras, J. Mathematical Physics
25(1984)2367-2373.


\bibitem{vergne1} A. Szenes and Mich\`{e}le Vergne, Residue formulae for vector partitions
and Euler-MacLaurin sums,Adv. in Appl. Math. 30(2003),
NO.1-2,295-342.

\bibitem{sturmfels} Bernd Sturmfels, On vector partition
functions, J. Combin. Theory Ser.A 72(2)(1995),302-309.



\bibitem{brion} Michel Brion and Mich\`{e}le Vergne, Residue
formulae, vector partition functions and lattice points in
rational polytopes, J. Amer. Math. Soc. 10(4)(1997),797-833.

\bibitem{xu}Zhiqiang Xu, Discrete Truncated Powers, Volume of Convex Polytopes and
Ehrhart Polynomials, http://arxiv.org/abs/math.CO/0505129.

\bibitem{popoviciu} Tiberiu Popoviciu, Asupra unei probleme de
patitie a numerelor, Acad. Republicii Populare Romane, Filiala
Cluj, Studii si cercetari stiintifice 4(1953),7-58.



\end{thebibliography}
\end{document}